\documentclass[reqno,openany]{article}
\newcommand{\beq}{\begin{equation}}
\newcommand{\eeq}{\end{equation}}
\newcommand{\barr}{\begin{eqnarray}}
\newcommand{\barn}{\begin{eqnarray*}}
\newcommand{\earr}{\end{eqnarray}}
\newcommand{\earn}{\end{eqnarray*}}

\begin{document}

\title{A statistical test for the Zipf's law by deviations from the Heaps' law
}
\author{Mikhail Chebunin\thanks{E-mail: chebuninmikhail@gmail.com, Sobolev Institute of Mathematics, 
Novosibirsk State University, Novosibirsk, Russia},
Artyom Kovalevskii \thanks{E-mail: artyom.kovalevskii@gmail.com, Novosibirsk State Technical University, 
Novosibirsk State University, Novosibirsk, Russia.  
The research was supported by RFBR grant 17-01-00683}}
\date{}
\maketitle

\begin{abstract}
We explore a probabilistic model of an artistic text: words of the text are chosen independently of each other
in accordance with a discrete probability distribution on an infinite dictionary. The words are enumerated 1, 2, $\ldots$,
and the probability of appearing the $i$'th word is asymptotically a power function. 
Bahadur proved that in this case the number of different words 
depends on the length of the text is asymptotically a power function, too. 
On the other hand, in the applied statistics community, there exist statements supported by empirical observations, the Zipf's and the Heaps' laws.
We highlight the links between Bahadur results and Zipf's/Heaps' laws, and 
introduce and analyse a corresponding statistical test.

\end{abstract}

\section{Introduction}

There is a countably infinite  dictionary where the words are numbered $1, \ 2, \ \ldots$. 
Words are chosen one-by-one independently of each other and accodingly to a discrete probability distribution
that is equivalent to a power law on the positive integers,
\begin{equation}\label{Zipf}
p_i\sim c i^{-1/\theta}, \ \ 0<\theta<1, \ \  c>0.
\end{equation}

Hereinafter for two positive sequences $\{a_j\}$ and $\{b_j\}$,
we write $a_j \sim b_j$ if $a_j/b_j \to 1$ as $j \to \infty$.

Denote by $ R_j $ the number of different words in the text of length $ j $.
Bahadur (1960) proved that
\begin{equation}\label{Bah}
{\bf E}R_{j} \sim c^{\theta} \Gamma(1-\theta) j^{\theta},
\end{equation}
where $\Gamma(\cdot)$ is the Euler gamma function.
Bahadur also proved
convergence in probability
$R_j/{\bf E}R_{j} \stackrel{p}{\to} 1$.

Karlin (1967) made next important steps. 
He proved that $R_j/{\bf E} R_{j} \stackrel{a.s.}{\to} 1$, which is equivalent to 
\begin{equation}\label{Heaps} 
R_{j} \sim C_1 j^{\theta} \ {\rm a.s.,}
\end{equation}
thanks to (\ref{Bah}).
Here $C_1=c^{\theta} \Gamma(1-\theta)$, coefficients $c$ and $\theta$ are from (\ref{Bah}), and 
$\Gamma(x)=\int_0^{\infty} y^{x-1} e^{-y}\,  dy$.

He proposed to consider texts of a random length using an independent rate-1 Poisson process $\Pi(t)$: the length of the text
grows in time and follows a Poisson distribution $\Pi(t)$ at time $t>0$.
Using this Poissanization, Karlin proved the Central Limit Theorem 
\begin{equation}\label{CLT} 
(R_j-{\bf E} R_{j})/\sqrt{{\bf Var} R_j} \Rightarrow N_{0,1}.
\end{equation}

Chebunin and Kovalevskii (2016) proved the Functional Central Limit Theorem, that is, the weak convergence of the process
$\{(R_{[nt]}-{\bf E} R_{[nt]})/\sqrt{{\bf Var} R_n}, \ 0\le t \le 1\}$
to a centered Gaussian
process $ Z $ with continuous a.s. sample paths and covariance function of the form
\[
K(s,t)=(s+t)^{\theta}-\max(s^{\theta}, t^{\theta}).
\]

Estimators of the parameter $\theta$ have been obtained by Nicholls (1978), Zakrevskaya and Kovalevskii (2001),
Guillou and Hall (2002), Ohannessian and Dahleh (2012), Chebunin (2014), Chebunin and Kovalevskii (2018).

On the other hand, relations (\ref{Zipf}) and (\ref{Heaps}) were observed empirically in the analysis of artistic texts.
Linguists call them the Zipf's law and the Heaps' law respectively.
The Zipf's law (Zipf, 1936) states the decrease in the frequencies of words depending on the rank 
in accordance with a power function.

The content of the Heaps' law  was initially proposed by Herdan (1960). This law was then popularized by Heaps (1978). 
The Heaps' law states that dependence of the number of different words from 
the text length is characterised by a power function.

Links between the Heaps' and the Zipf's laws have been studied (empirically and in other contexts) 
by van Leijenhorst and van der Weide (2005),
Serrano et al. (2009), Bernhardsson et al. (2009), Eliazar (2011), Baeza-Yates $\&$ Navarro (2013), ets.

As far as we are aware of, statistical test of fit for Zipf's law is still the problem. 
We are unaware of any mathematically correct statistical goodness-of-fit test Zipf's law.  
Altmann and Gerlach (2016) emphasize incorrectness of a number of statistical tests proposed earlier.

Note that analysis of very long texts and text sequences shows significant deviations
from the Zipf model (Petersen et al., 2012). Gerlach and Altmann (2013) proposed a modified model for explanations
of these deviations.

The present paper proposes a new theoretically supported test for the Zipf's law. We introduce a new class of estimates 
that is based on the sequence
$ (R_1, \ldots, R_n) $. We form an empirical process and prove its weak convergence to a centered Gaussian process. 
We calculate the covariance function of this limiting process. Then we construct a test
of the omega-squared type. Calculation of the limiting
distribution of the test statistics is based on the corresponding classical result of Smirnov (1937).
One can calculate the distribution using results of Deheuvels and Martynov (1996).

The rest of the paper is organised as follows. We introduce
an estimator for parameter $\theta$ in Section 2, a test for known $\theta$ in Section 3, and 
a test for unknown $\theta$ (that is, 
the test for the Zipf's law) in Section 4. Proofs are given in Section 5.

\section{Parameter's estimation }

From (\ref{Heaps}),  we have $\log R_j\sim{\theta}\log j$  a.s.
Therefore, we may propose the following estimator for parameter $\theta$:
\[
\widehat{\theta}=\int_0^1  \log^+ R_{[nt]} \, d A(t)
\]
with function $A(\cdot)$ such that
\begin{equation}\label{intalog}
\int_0^1 \log t \, dA(t) =1, \ \ A(0)=A(+0)=A(1) =0,
\end{equation}
here $\log^+ x=\max(\log x, \, 0)$.
We assume $A(\cdot)$ to be the sum of a step function and a piecewise continuosly differentiable function on $[0,\, 1]$.

{\bf Theorem 2.1} 
If $p_i =i^{-1/\theta} l(i,\theta)$, $\theta\in [0,1]$, and if $l(x,\theta)$ 
is a slowly varying function as $x \to \infty$, then 
the estimator $\widehat{\theta}$ is strongly consistent.

We need extra conditions to obtain the asymptotic normality of $\widehat{\theta}$.
 
{\bf Theorem 2.2}
Let 
 $A(t)=0$, $t \in [0,\, \delta]$ for some $\delta\in(0,\, 1)$,  and $p_i=ci^{-1/\theta}(1+o(i^{-1/2}))$,
 $\theta \in (0,1)$. Then
\[
\sqrt{{\bf E}R_n} (\widehat{\theta}-\theta) - \int_0^1  t^{-\theta} Z_n(t) \, dA(t) \to_p 0.
\]

From Theorem 2.2, it follows that  $\widehat{\theta}$ converges to $\theta$ with rate $({\bf E}R_n)^{-1/2}$, and 
normal random variable $ \int_0^1  t^{-\theta} Z_{\theta}(t) \, dA(t)$ has variance
$ \int_0^1 \int_0^1 (st)^{-\theta} K(s,t) \, dA(s) \, dA(t)$.

{\bf Example 2.1} Take
\[
A(t) = \left\{
\begin{array}{ll}
0, & 0 \le t \le 1/2; \\
-(\log 2)^{-1}, & 1/2<t<1; \\
0, & t=1.
\end{array}
\right.
\]

Then
\[
\widehat{\theta} = \log_2 (R_n/R_{[n/2]}), \ \ n  \ge 2.
\]

Note that, in this example, for any function $g$ on $[0,1]$, 
\[
\int_0^1 g(t)\, dA(t)=\frac{g(1)-g(1/2)}{\log 2}.
\]

\section{Test for a known rate}

Let $0<\theta<1$ be known. We introduce an {\em empirical bridge} $Z_n^0$ (Kovalevskii and Shatalin, 2015) as follows.
\[
Z_n^0(k/n)=\left(R_k-(k/n)^{\theta}R_n\right)/\sqrt{R_n},
\]
$0\leq k \leq n$, where $R_0=0$. We construct a piecewise linear approximation:
for $0\leq t \leq 1/n$ and $0\leq k \leq n-1$,
\[
Z_n^0\left(\frac{k}{n}+t\right)=Z_n^0(k/n)+nt\left(Z_n^0((k+1)/n) - Z_n^0(k/n)\right).
\] 

{\bf Theorem 3.1} Under the assumptions of Theorem 2.2,
\[
\sup_{0\le t \le 1} |Z_n^0(t)-(Z_n(t)-t^{\theta}Z_n(1))| \to 0 \ \mbox{a.s.}
\]

Let C(0,1) be the set of all continious functions on $[0,1]$ with the uniform metric
$\rho(x,y)=\max_{t\in[0,1]}|x(t)-y(t)|$.
By the FCLT of Chebunin $\&$ Kovalevskii (2016), we have

{\bf Corollary 3.1} Under the assumptions of Theorem 2.2,
$Z_n^0$ converges weakly in $C(0,1)$ to a Gaussian process $Z_{\theta}^0$ 
that can be represented as $Z_{\theta}^0(t)=Z_{\theta}(t) - t^{\theta} Z_{\theta}(1)$, $0 \leq t \leq 1$. 
Its correlation function is given by
\[
K^0(s,t)={\bf E} Z_{\theta}^0(s) Z_{\theta}^0(t)=K(s,t)-s^{\theta}K(1,t)-t^{\theta}K(s,1)+s^{\theta}t^{\theta}K(1,1).
\]

Now we show how to implement the goodness-of-fit test
in this case.
Let $W_n^2=\int\limits_0^1 \left(Z_n^0(t)\right)^2 dt$. It is equal to 
\begin{equation}\label{wn2}
W_n^2=\frac{1}{3n}\sum_{k=1}^{n-1} Z_n^0\left(\frac{k}{n}\right)\left( 2Z_n^0\left(\frac{k}{n}\right)+Z_n^0\left(\frac{k+1}{n}\right) \right).
\end{equation}

Then $W_n^2$ converges weakly to $W_{\theta}^2= \int\limits_0^1 \left(Z_{\theta}^0(t)\right)^2 dt$.

So the test rejects the basic hypothesis if $W_n^2 \geq C$. The p-value of the test is $1-F_{\theta}(W_{n,obs}^2)$.
Here $F_{\theta}$ is the cumulative distribution function of $W_{\theta}^2$ and $W_{n,obs}^2$ is a concrete value of $W_{n}^2$
 for observations under consideration.

One can estimate $F_{\theta}$ by simulations or find it explicitely using the  Smirnov's formula (Smirnov, 1937): 
if $W_{\theta}^2=\sum_{k=1}^{\infty}\frac{\eta^2_k}{\lambda_k}$,  
$\eta_1,\eta_2,\ldots$ are independent and have standard normal distribution, $0<\lambda_1<\lambda_2<\ldots$, then 
\begin{equation}\label{cdf}
F_{\theta}(x)=1+\frac{1}{\pi} \sum_{k=1}^{\infty} (-1)^k \int_{\lambda_{2k-1}}^{\lambda_{2k}}\frac{e^{-\lambda x/2}}{\sqrt{-D(\lambda)}}
\cdot \frac{d\lambda}{\lambda}, \ x>0,
\end{equation}
\[
D(\lambda)=\prod_{k=1}^{\infty} \left( 1- \frac{\lambda}{\lambda_k}\right).
\]

The integrals in the RHS of (\ref{cdf}) must tend to 0 monotonically as $k \to \infty$, 
and $\lambda_k^{-1}$ are the eigenvalues of kernel $K^0$
(see Martynov (1973), Chapter 3).

\section{Test for an  unknown rate}

Let us introduce the process $\widehat{Z}_n$:
\[
\widehat{Z}_n(k/n)=\left(R_k-(k/n)^{\widehat{\theta}}R_n\right)/\sqrt{R_n},
\]
$0\leq k \leq n$.
As for ${Z}_n^0$, let
for $0\leq t \leq 1/n$ and $0\leq k \leq n-1$
\[
\widehat{Z}_n\left(\frac{k}{n}+t\right)=\widehat{Z}_n(k/n)+nt\left(\widehat{Z}_n((k+1)/n) - \widehat{Z}_n(k/n)\right).
\]

{\bf Theorem 4.1} Under assumptions of Theorem 2.2, 
$\widehat{Z}_n$ converge weakly to $\widehat{Z}_{\theta}$ as $n \to \infty$, where 
\[
\widehat{Z}_{\theta}(t) ={Z}_{\theta}^0(t) - t^{\theta}\log  t 
\int_0^1 u^{-\theta} Z_{\theta}(u)\,dA(u).
\]

{\bf Corollary 4.1}  Assume the conditions of Theorem 2.2 to hold. 
Let $\widehat{W}_n^2=\int\limits_0^1 \left(\widehat{Z}_n(t)\right)^2 dt$.
Then $\widehat{W}_n^2$ converges weakly to $\widehat{W}_{\theta}^2= \int\limits_0^1 \left(\widehat{Z}_{\theta}(t)\right)^2 dt$.

Similarly to (\ref{wn2}), $\widehat{W}_n^2$ has the following representation 
\[
\widehat{W}_n^2=\frac{1}{3n}\sum_{k=1}^{n-1} \widehat{Z}_n\left(\frac{k}{n}\right)\left( 2\widehat{Z}_n\left(\frac{k}{n}\right)+
\widehat{Z}_n\left(\frac{k+1}{n}\right) \right).
\]

The p-value of the goodness-of fit test is $1-\widehat{F}_{\theta}(\widehat{W}_{n,obs}^2)$.
Here $\widehat{F}_{\theta}$ is the cumulative distribution function of $\widehat{W}_{\theta}^2$, and $\widehat{W}_{n,obs}^2$ 
is the observed value of $\widehat{W}_n^2$. Further, the function
$\widehat{F}_{\theta}$ can be found using the approach from Section 3, with replacing 
$\lambda_k$ by $\widehat{\lambda}_k$ 
in the Smirnov's formula, and
$\widehat{\lambda}_k$ are the eigenvalues of the kernel $\widehat{K}(s,t)={\bf E} \widehat{Z}_{\theta}(s) \widehat{Z}_{\theta}(t) $.

\section{Proofs}

{\em Proof of Theorem 2.1}

Let
\[
\alpha(x)=\max\{j|\ p_j\geq 1/x\}.
\]

Since $p_i i^{1/\theta}$ is a slowly varying function as $i\to \infty$, we have 
$\alpha(x)=x^{\theta} L(x,\theta)$, 
$L(x,\theta)$ is a slowly varying function as $x \to \infty$ (Karlin, 1967).

Take a sequence $\{\delta_n\}$ such that 
$n \delta_n \to \infty$, $delta_n\log n\to 0$ (for example, $\delta_n=\log ^{-2} n$).  Then for sufficiently large $n$
\[
\left|\int_0^{\delta_n} \log  R_{[nt]} \, dA(t)\right|\le_{a.s.} 
\max\limits_{0\le t\le \delta_n} |A'(t)| \delta_n \log  n\to 0.
\]

The rest of the integral is
\[
\int_{\delta_n}^1 \log  R_{[nt]} \, dA(t) =\int_{\delta_n}^1  \log  \frac{R_{[nt]}}{{\bf E}R_{[nt]}} \, dA(t)
+ \int_{\delta_n}^1  \log  {\bf E}R_{[nt]} \, dA(t).
\]

We prove the a.s. convergence of the first integral in RHS to 0 a.s., and then convergence of the second one to $\theta$.

By the SLLN,  $\log (R_j/{\bf E}R_j) \to 0$ a.s. as $j\to \infty$.
Therefore, for any $\varepsilon >0$,
\[
\lim_{n \to \infty}  
{\bf P}\left(\sup_{j \ge n \delta_n}\left|\log  \left(\frac{R_j}{{\bf E}R_j}\right)\right|\ge \varepsilon\right)=0.
\]

Therefore
\[
 {\bf P}\left( \sup_{k \ge n }\left|\int_{\delta_n}^1  \log  \frac{R_{[kt]}}{{\bf E}R_{[kt]}}\, dA(t) \right|\ge 
\varepsilon \right)\le
{\bf P}\left( \int_{\delta_n}^1 \sup_{k \ge n } \left| \log  \frac{R_{[kt]}}{{\bf E}R_{[kt]}}\right|\, |dA(t)| \ge 
\varepsilon \right)
\]
\[
 ={\bf P}\left( \sup_{k \ge n\delta_n }\left| \log  \frac{R_{k}}{{\bf E}R_{k}} \right|\ge 
\frac{\varepsilon}{ \int_{\delta_n}^1 |dA(t)|}\right)\to 0 \ \ \textrm{as} \ \ n\to\infty,
\]
that is, $\int_{\delta_n}^1 \log  \frac{R_{[nt]}}{{\bf E}R_{[nt]}} \, dA(t) \to 0$ a.s.

We have ${\bf E} R_j=j^{\theta} L(j)$, where $L(j)$ is a slowly varying function (Karlin, 1967).
So, uniformly on $t \ge \delta_n>0$,
\[
\frac{{\bf E} R_{[nt]}}{ (nt)^{\theta}L(nt)} \to 1,
\]
\[
\log  \frac{{\bf E} R_{[nt]}}{ (nt)^{\theta}L(nt)}  \to 0.
\]

Therefore
\[
\int_{\delta_n}^1  \log  \frac{{\bf E} R_{[nt]}}{ (nt)^{\theta}L(nt)} \, dA(t) \to 0.
\]

However, $L(nt)=(nt)^{o(1)}$ for $t\ge \delta_n$, so from (\ref{intalog})   
we have 
\[
\int_{\delta_n}^1  \log  ( (nt)^{\theta}L(nt)) \, dA(t) = \int_{\delta_n}^1  \log  ( (nt)^{\theta+o(1)}) \, dA(t) = 
(\theta+o(1))\int_{\delta_n}^1  \log  (nt) \, dA(t) 
\]
\[
 = 
(\theta+o(1))\int_{0}^1 \log  (nt) \, dA(t) - (\theta+o(1))\int_{0}^{\delta_n}  \log  (nt) \, dA(t) = \theta+o(1).
\]

Then
\[
\int_{\delta_n}^1  \log  {\bf E} R_{[nt]} \, dA(t) \to \theta,
\]
and $\widehat{\theta}\to \theta$ a.s. The proof is complete.

{\em Proof of Theorem 2.2}

Since 
$p_i=ci^{-1/\theta}(1+o(i^{-1/2}))$, we have 
\begin{equation}\label{ern}
{\bf E}R_n=C_1 n^\theta + o(n^\frac\theta2)
\end{equation}
(Chebinin and Kovalevskii, 2017, Lemma 1, 2). Recall that
\[
Z_n(t)=\frac{R_{[nt]} - {\bf E} R_{[nt]}}{\sqrt{{\bf E}R_n}}.
\]

Let
\[
Z^*_n(t)=\frac{R_{[nt]} - {\bf E} R_{[nt]}}{{\bf E}R_{[nt]}}.
\]

Then
\[
\sqrt{{\bf E}R_n} \left( \int_0^1  \log  R_{[nt]} \, dA(t) -\theta \right) - \int_0^1 t^{-\theta} Z_n(t) \, dA(t) 
\]
\[
=\sqrt{{\bf E}R_n} \left(\int_0^1 \left(\log  C_1 (nt)^\theta+\log \frac{ {\bf E}R_{[nt]}}{C_1 (nt)^\theta}  
+\log  (1+Z_n^*(t))-t^{-\theta} \frac{Z_n(t)}{\sqrt{{\bf E}R_n} }\right)\, dA(t) -\theta\right) 
\]
\[
=\sqrt{{\bf E}R_n} \int_0^1  \left(\log \frac{ {\bf E}R_{[nt]}}{C_1 (nt)^\theta}  
+\log  (1+Z_n^*(t))-Z_n^*(t)+Z_n^*(t)-t^{-\theta} \frac{Z_n(t)}{\sqrt{{\bf E}R_n}}\right)\, dA(t). 
\]

For any $t\in[\delta,1]$, $Z_n^*(t)\to 0$ a.s. , so $\log  (1+Z_n^*(t))-Z_n^*(t)\sim - (Z_n^*(t))^2/2 $ a.s.
We have ${\bf E}R_n= n^\theta L(n)$, so, for any $t\in[\delta,1]$,
\[ 
\frac{({\bf E}R_n)^\frac32}{({\bf E}R_{[nt]})^2}
= \frac{(n^\theta L(n))^\frac32}{((nt)^\theta L(nt))^2}=O\left(\frac1{n^\frac{\theta}{2} \sqrt{ L(n)}}\right), 
\ \ \ \frac{t^\theta {\bf E}R_n}{{\bf E}R_{[nt]}} =\frac{L(n)}{L(nt)}=1+o(1).
\]

Note that
\[
\sqrt{{\bf E}R_n} \int_0^1  (Z_n^*(t))^2\, dA(t)=  \int_0^1  \frac{({\bf E}R_n)^\frac32}{({\bf E}R_{[nt]})^2} (Z_n(t))^2\, dA(t)
\]
\[
=o(n^{-\theta/4}) \int_0^1 (Z_n(t))^2\, dA(t)\to_p 0,
\]
so
\[
\sqrt{{\bf E}R_n} \int_0^1 \left(Z_n^*(t)-t^{-\theta} \frac{Z_n(t)}{\sqrt{{\bf E}R_n}}\right)\, dA(t)=  
 \int_0^1  \left(\frac{t^\theta {\bf E}R_n}{{\bf E}R_{[nt]}}-1\right)t^{-\theta} Z_n(t)\, dA(t)
\]
\[
=  
o(1)\int_0^1  t^{-\theta} Z_n(t)\, dA(t)\to 0.
\]

Since ${\bf E}R_n=C_1 n^\theta + o(n^\frac\theta2)$, we have
\[
\sqrt{{\bf E}R_n} \int_0^1  \log \frac{ {\bf E}R_{[nt]}}{C_1 (nt)^\theta} \, dA(t)
=\sqrt{{\bf E}R_n} \int_0^1  \log (1+o(n^{-\frac\theta2})) \, dA(t) =o(1)\int_0^1 dA(t)\to 0. 
\]

The proof is complete.

{\it Proof of Theorem 3.1}

Due to inequality
\[
\sup_{0\leq t \leq 1} \left|Z_n^0(t) - \frac{R_{[nt]}-t^{\theta}R_n}{\sqrt{R_n}}\right|\leq \frac{1}{n},
\]
convergence (\ref{Heaps}) it is enough to prove that
\[
\sup_{0\leq t \leq 1} \left|\frac{R_{[nt]}-t^{\theta}R_n}{\sqrt{{\bf E}R_n}} - (Z_n(t)-t^{\theta}Z_n(t))\right|\to 0.
\]
The latter follows from the definition of $Z_n$ and (\ref{ern}).
The proof is complete.

{\it Proof of Theorem 4.1}

By the definition,
\[
\widehat{Z}_n(k/n)={Z}_n^0(k/n)-\sqrt{R_n}\left((k/n)^{\widehat{\theta}}-(k/n)^{{\theta}}\right),
\] 

If $k \le nt \le k+1$, then (see the proof of Theorem 3.1)
\[
\frac{R_k-\left(\frac{k+1}{n}\right)^{\widehat{\theta}} R_n}{\sqrt{R_n}}\le \widehat{Z}_n(t) \le 
\frac{R_{k+1}-\left(\frac{k}{n}\right)^{\widehat{\theta}} R_n}{\sqrt{R_n}}.
\]

So
\[
\sup_{0\leq t \leq 1} \left|\widehat{Z}_n(t)-{Z}_n^0(t)+\sqrt{R_n}\left(t^{\widehat{\theta}}-t^{{\theta}}\right)\right|
\]
\[
\le \frac{\sqrt{R_n}}{n^{\theta}} +\frac{\sqrt{R_n}}{n^{\widehat{\theta}}}+ \frac{1}{\sqrt{R_n}} \to 0 \ \mbox{a.s.},
\]
and
\[
\sup_{0\leq t \leq 1} \left| \widehat{Z}_n(t)-{Z}_n^0(t)+\sqrt{{\bf E} R_n}\left(t^{\widehat{\theta}}-t^{{\theta}}\right)\right|
 \to 0 \ \mbox{a.s.}
\]

Since $\widehat{\theta}\to \theta$ a.s., we have 
$
t^{\widehat{\theta}}-t^{{\theta}}\sim (\widehat{\theta}-\theta)t^{\theta}\log  t 
$
a.s. for any  $0<t\le 1$. 

From Theorem 2.2 and Theorem 3.1, we have 
\[
\widehat{Z}_{\theta}(t) ={Z}_{\theta}^0(t) -t^{\theta}\log  t 
\int_0^1 u^{-\theta} Z_{\theta}(u)\,dA(u).
\]

The proof is complete. 

\bigskip

{\bf Acknowledgements}

The research was supported by RFBR grant 17-01-00683.
The authors like to thank Sergey Foss 
for helpful and constructive comments and suggestions.

\bigskip


\bigskip

\footnotesize

{\sc Altmann, E. G., and Gerlach, M.,}, 2016.
Statistical laws in linguistics. In: M. Degli Esposti et al. (eds.), Creativity and Universality in Language,
Lecture Notes in Morphogenesis.

{\sc Baeza-Yates, R., and Navarro, G.}, 2000.  Block Addressing Indices
for Approximate Text Retrieval, J. Am. Soc. Inf. Sci. 51, 69.

{\sc Bahadur, R. R.}, 1960. On the number of distinct values in a large sample from an infinite discrete distribution.
Proceedings of the National Institute of Sciences of India, 26A, Supp. II, 67--75.

{\sc Bernhardsson, S., Correa da Rocha, L. E., Minnhagen, P.}, 2009. The Meta Book and Size-Dependent
Properties of Written Language, New J. Phys. 11, 123015.

{\sc Chebunin, M. G.}, 2014. Estimation of parameters of probabilistic models which is based on the number of different elements 
in a sample. Sib. Zh. Ind. Mat., 17:3, 135--147 (in Russian).

{\sc Chebunin, M., Kovalevskii, A.}, 2016. 
Functional central limit theorems for certain statistics in an infinite urn scheme. 
Statistics and Probability Letters, V. 119, 344--348.

{\sc Chebunin, M., Kovalevskii, A.}, 2018.
Asymptotically normal estimators for Zipf's law. Sankhya A.

{\sc Deheuvels, P., and Martynov, G. V.}, 1996. 
Cramer-von mises-type tests with applications to tests of independence for multivariate extreme-value distributions.
Communications in Statistics --- Theory and Methods, 25:4, 871--908.

{\sc Gerlach, M., and  Altmann, E. G.}, 2013.
Stochastic Model for the Vocabulary Growth in Natural Languages. 
Physical Review X 3, 021006.

{\sc Guillou, A., Hall, P.}, 2002.
A diagnostic for selecting the threshold in extreme value analysis.
Journal of the Royal Statistical Society: Series B. Volume 63, Issue 2, 293--305.



{\sc Heaps, H. S.}, 1978.
Information Retrieval: Computational and Theoretical Aspects, Academic Press.

{\sc Herdan, G.}, 1960.
 Type-token mathematics, The Hague: Mouton.

{\sc Eliazar, I.}, 2011. The Growth Statistics of Zipfian Ensembles:
Beyond Heaps' Law, Physica (Amsterdam) 390, 3189.


{\sc Karlin, S.}, 1967. Central Limit Theorems for Certain Infinite Urn Schemes. 
Jounal of Mathematics and Mechanics, Vol. 17, No. 4,  373--401.

{\sc Kovalevskii, A. P.,  Shatalin, E. V.}, 2015. Asymptotics of sums of residuals of one-parameter linear regression 
on order statistics, Theory of probability and its applications, Vol. 59, No. 3, 375--387.

{\sc van Leijenhorst, D. C., van der Weide, T. P.}, 2005.  A Formal
Derivation of Heaps' Law, Information Sciences (NY)
170, 263.

{\sc Martynov, G. V.}, 1978. The omega squared test. Nauka, Moscow (in Russian).


{\sc Nicholls, P. T.}, 1987.  
Estimation of Zipf parameters. J. Am. Soc. Inf. Sci., V. 38, 443--445.

{\sc Ohannessian, M. I., Dahleh, M. A.}, 2012.
Rare probability estimation under regularly varying heavy tails,   
Proceedings of the 25th Annual Conference on Learning Theory, PMLR 23:21.1--21.24. 

{\sc Petersen, A. M., Tenenbaum,  J. N.,  Havlin, S., Stanley,  H. E., Perc,  M.}, 2012 .
Languages cool as they expand: Allometric scaling and the decreasing need for new words.
Scientific Reports 2, Article No. 943.

{\sc Serrano, M. A., Flammini, A., Menczer, F.}, 2009. Modeling
Statistical Properties of Written Text, PLoS ONE 4, e5372.

{ \sc Smirnov, N.V.}, 1937. On the omega-squared distribution, Mat. Sb.  2,  973--993 (in Russian).

{\sc Zakrevskaya, N. S.,  Kovalevskii, A. P.}, 2001. One-parameter probabilistic models of text statistics. 
Sib. Zh. Ind. Mat., 4:2, 142--153 (in Russian).

{\sc Zipf, G. K.}, 1936. The Psycho-Biology of Language. Routledge, London.

\end{document}